\documentclass[12pt]{amsart}
\usepackage{geometry} % see geometry.pdf on how to lay out the page. There's lots.
\usepackage{amsfonts}
\usepackage{amssymb}
\usepackage{graphicx}
\title{Some error estimates for the DEC method in the plane}
\author{Ruben Carrillo \and Miguel Angel Moreles \and Rafael Herrera}
\address{R. Carrillo, M. A. Moreles, Rafael Herrera \\
                 \email{moreles@cimat.mx}\\
              Centro de Investigaci\'{o}n en Matem\'{a}ticas\\ Jalisco s/n, Valenciana\\
Guanajuato, GTO 36240,Mexico           
}
\date{} % delete this line to display the current date

\begin{document}
\begin{abstract}
We show that the Discrete Exterior Calculus (DEC) method can be cast as the earlier box method for the Poisson problem in the plane. Consequently, error estimates are established, proving that the DEC method is comparable to the Finite Element Method with linear elements. We also discuss some virtues, others than convergence, of the DEC method.
\end{abstract}
\maketitle
\tableofcontents

\section{Introduction}
Discrete Exterior Calculus is introduced in \cite{HiraniThesis} as a discrete version of Exterior Differential Calculus. Most mathematical objects in the continuous theory and the potential applications, have been explored actively with DEC. An application of interest in this work, is the numerical approximation of Partial Differential Equations (PDE). The attractive feature is that differential forms are the natural language for physical laws modeled by PDE. 

At present, the method has proven successful. For instance, in \cite{Hirani_K_N} it is applied to solve Darcy flow and Poisson's equation. As illustrated in \cite{Mohamed_etal_16}, the method is natural for curved meshes. Therein, they show an application to the incompressible Navier-Stokes equations over surface simplicial meshes.

These numerical results are proof that the DEC method is a welcomed addition to PDE numerics. But, to our knowledge, theoretical results on convergence  are scarce at best. There are works exploring numerical convergence such as  the recent by Mohamed et al \cite{Mohamed_etal_18}.

In \cite{griebel_etal}, the authors show that in simple cases (e.g. flat geometry and regular meshes), the equations
resulting from DEC are equivalent to classical numerical methods, finite difference or finite volume. It is well known that the latter involves a dual mesh in close analogy to the DEC method. This leads to the earlier paper of \cite{BankRose}, where the box method, a finite volume type method, is analyzed. It is straightforward  to join the dots and see that the box and DEC methods coincide for the Poisson equation and other simple self
adjoint variants. Our purpose is to provide a proof of this fact while introducing the basics of the DEC method. Consequently a convergence result is established, namely the DEC method converges comparable to the Finite Element Method with linear elements. We do not aim for generality, thus our description of the DEC and BOX methods just provides enough details as to make the connection.

The outline is as follows.

In Section 2 we review the BOX method following closely \cite{BankRose}. Local aspects, box-wise, are stressed. In Section 3 we introduce cell complexes in DEC theory, and interpret the meshes in the FEM and BOX methods as the primal and dual cell complexes respectively.
The discrete operators to approximate the Poisson equation are presented in Section 4. The proof of DEC equals BOX method is also carried out. Next, error estimates and first order convergence are established. A case for DEC is made, when classical methods are
limited by geometric issues. 

\section{The box method for the Poisson equation in 2D}

This section is a summary of \cite{BankRose}. 

\subsection{Meshes}

Let $\Omega$ be a bounded polygonal domain in $\mathbb{R}^2$ and $\partial\Omega$ its boundary. Let $\mathcal{T}$ a Finte Element triangulation with vertices $x_i,\, i=1,2,\ldots,n.$

To each vertex $x_i$ associate a region $\Omega_i$ consisting of those triangles $t\in\mathcal{T}$ with $x_i$ as a vertex, as in Figure 1.

\begin{figure}[htbp] %  figure placement: here, top, bottom, or page
   \centering
   \includegraphics[width=2in]{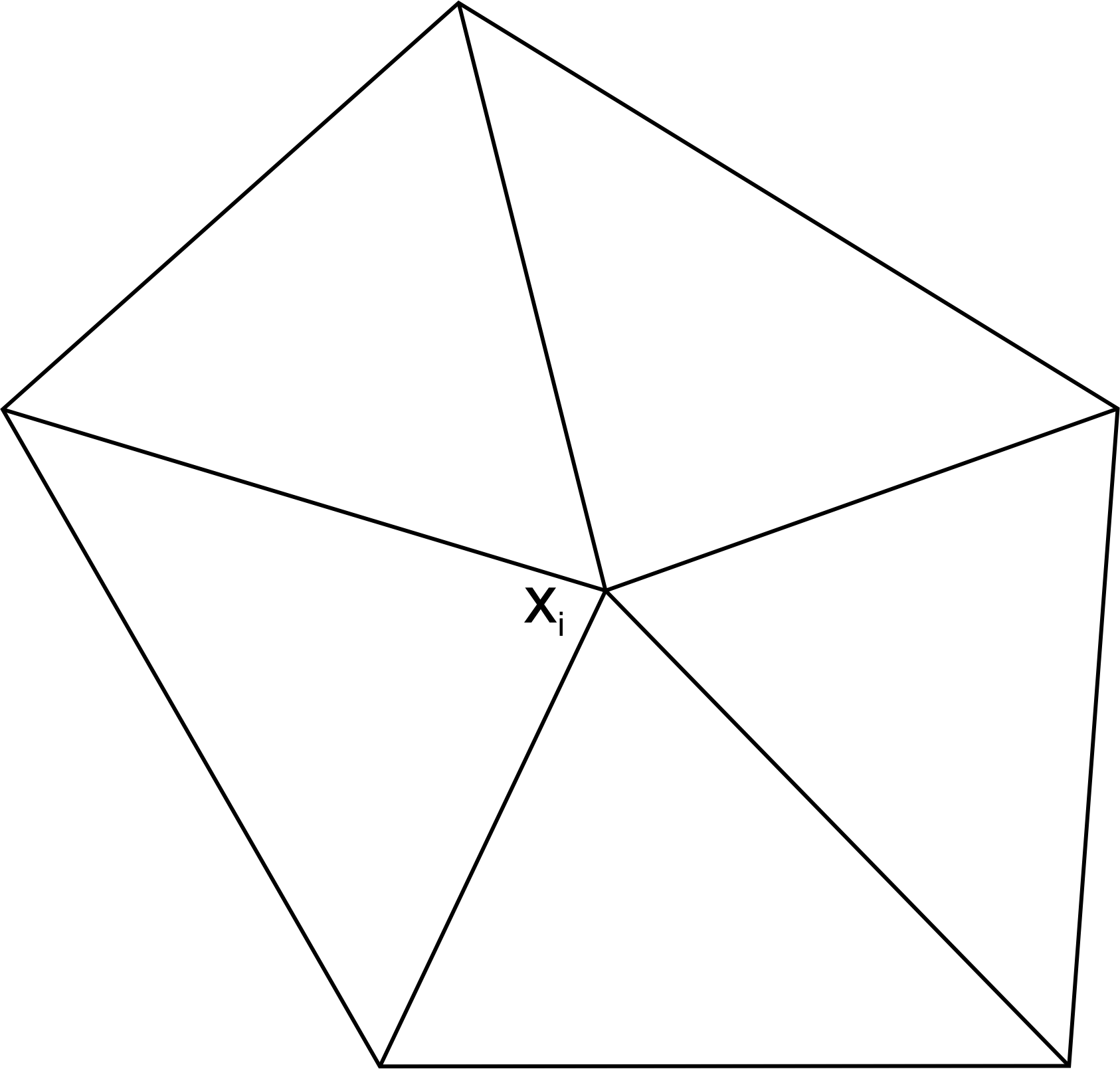} 
   \caption{Region $\Omega_i$.}
   \label{fig:example}
\end{figure}

A box is a basic element of a dual mesh constructed as follows, see Figure 2. For each triangle $t\in\mathcal{T}$ select a point $C\in\bar{t}$, later we shall restrict the mesh in order to select the circumcenter.

\begin{figure}[htbp] %  figure placement: here, top, bottom, or page
   \centering
   \includegraphics[width=2in]{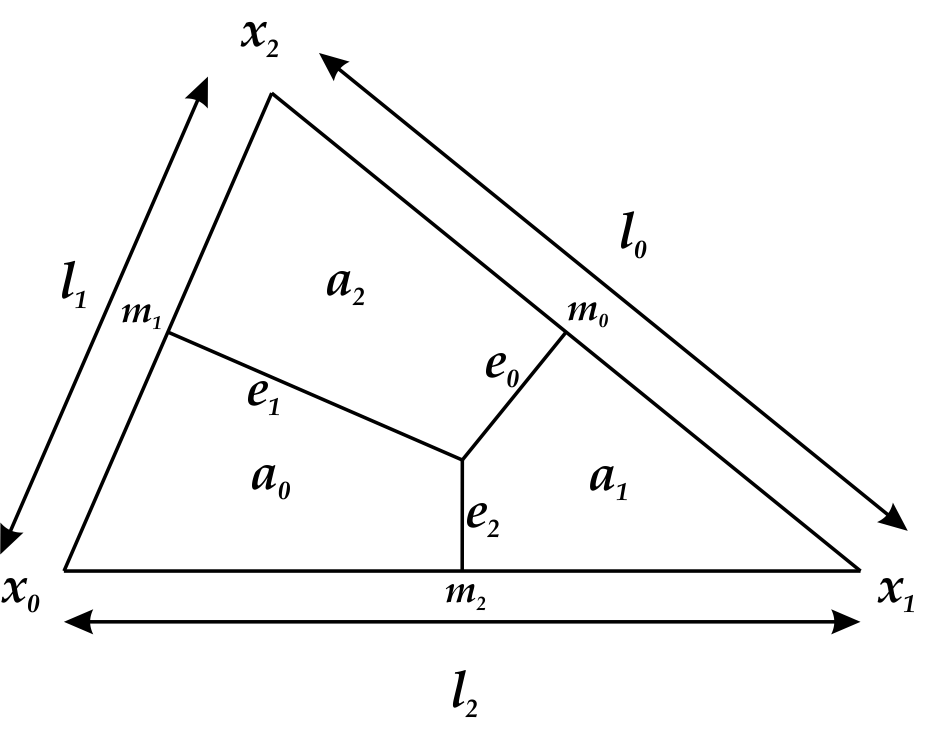} 
   \caption{Triangle $x_0x_1x_2$ and box elements}
   \label{fig:example}
\end{figure}

The point $C$ is connected with straight line segments ($e_0$, $e_1$, $e_2$) to the edge midpoints of $t$, ($m_0$, $m_1$, $m_2$). We obtain a partition of $t$ in three subregions ($a_0$, $a_1$, $a_2$). Sometimes for simplicity, $e_i$ $(a_i)$ denote both the set and its length (area).

We are led to a dual mesh $\mathcal{B}$ for $\mathcal{T}$ made out of boxes. For each vertex $x_i$ there is a corresponding box $b_i\subset\Omega_i$, consisting of the union of subregions in $\Omega_i$ which have $x_i$ as a corner, see Figure 3.

\begin{figure}[htbp] %  figure placement: here, top, bottom, or page
   \centering
   \includegraphics[width=2in]{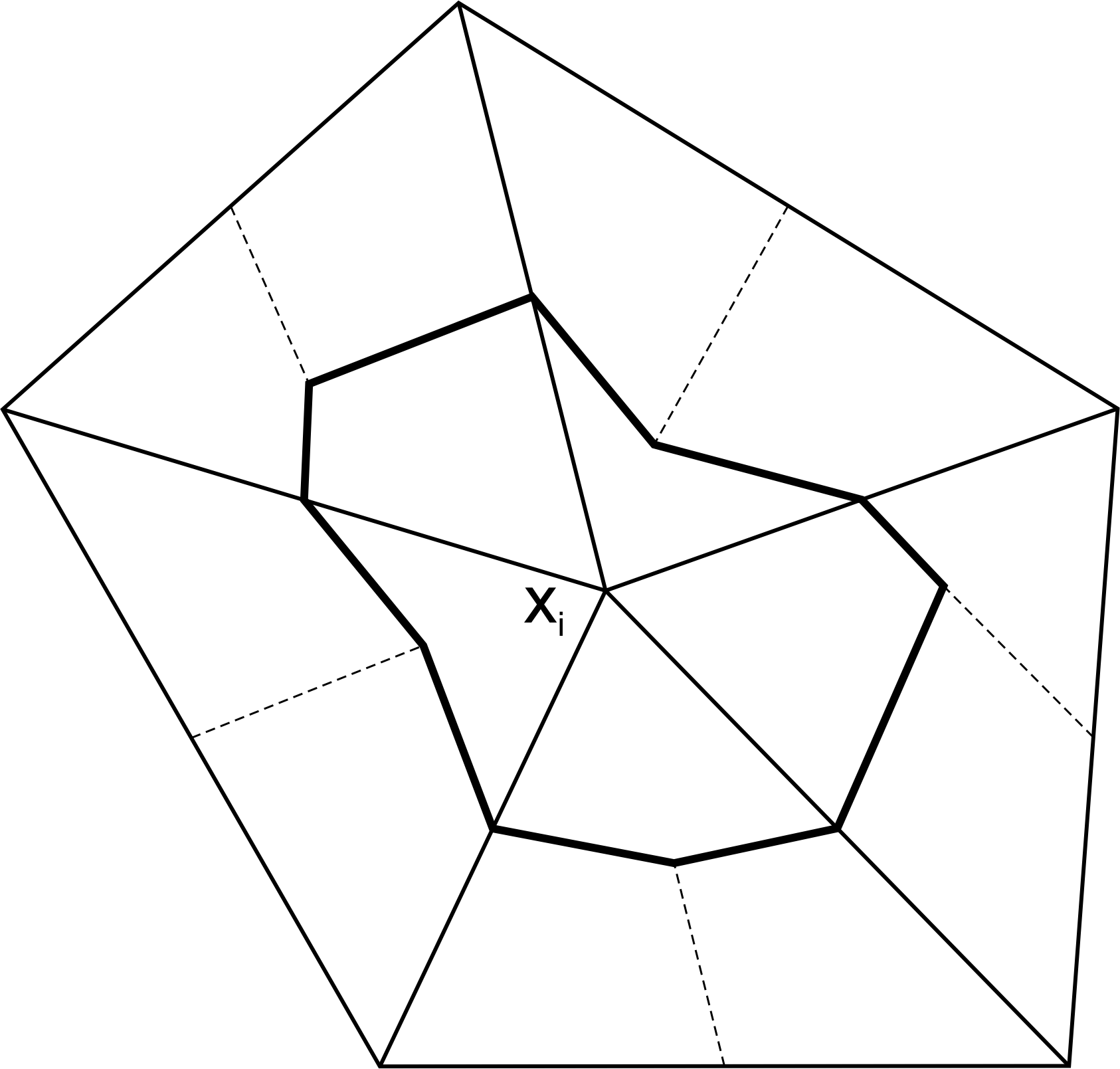} 
   \caption{Box $b_i$}
   \label{fig:example}
\end{figure}

The straight line segments from the selected points $C$ to the edge midpoints are collected in a set of edges, $\mathcal{E}$. 

We stress the following.

\bigskip

\noindent\textbf{(2.1) Definition. } For a box $b\in\mathcal{B}$, let us define its boundary $\eth b$ as the intersection of the topological boundary with the set of edges. 

\bigskip

Notice that only for vertices in the open set $\Omega$ the boundary is the topological one.

\subsection{Function spaces}

Let $L^2(R)$ and $H^1(R)$ denote the usual Sobolev spaces equipped with the norms
\[
\Vert u\Vert^2_{L^2(R)}=\int_R\vert u\vert^2dx,\quad 
\Vert u\Vert^2_{H^1(R)}=\Vert u\Vert^2_{L^2(R)}+\Vert \nabla u\Vert^2_{L^2(R)},
\]
where 
\[
\Vert \nabla u\Vert^2_{L^2(G)}=\int_G\vert\nabla u\cdot\nabla u\vert dx.
\]

If $u\in H^1_0(R)$ the $H^1$ norm is equivalent to the \emph{energy norm}
\[
\Vert u\Vert_E=\Vert \nabla u\Vert^2_{L^2(G)}.
\]
In $H^1_0(R)$ we also define bilinear form $a(u,v)$
\[
a(u,v)=\int_\Omega\nabla u\cdot\nabla vdx.
\]

Following the linear Finite Element Method, let $P^1(\mathcal{T})\subset H^1(\Omega)$ be the subspace of continuous piecewise linear polynomials associated with $\mathcal{T}$. Also consider
\[
S^0(\mathcal{B})=\lbrace v\in L^2(\Omega)\,:\, v\in H^1(b), b\in\mathcal{B} \rbrace.
\]
The subspace $P^0(\mathcal{B})$ of $S^0(\mathcal{B})$, is the space of discontinuous piecewise constants with respect to the boxes.

Let $\lbrace \phi_i\rbrace$ the usual nodal basis for $P^1(\mathcal{T})$ satisfying
\[
\phi_i(x_j)=\delta_{ij}.
\]

Let $\lbrace \bar{\phi}_i\rbrace$ the basis for $P^0(\mathcal{B})$ consisting of characteristic functions for $b_i$.

If $u\in P^1(\mathcal{T})$, then for some unique scalars $a_i$ we have
\[
u(x)=\sum_ia_i\phi(x).
\]

Easily, the map
\[
G:P^1(\mathcal{T})\to P^0(\mathcal{B}),
\]
given by
\[
G(u)(x)=\bar{u}(x)=\sum_ia_i\bar{\phi}_i(x),
\]
is bijective. Moreover, $u$ and $\bar{u}$ take the same value at the vertices of $\mathcal{T}$.

\bigskip

Let $P^0_0(\mathcal{B})$ be the subspace of $P^0(\mathcal{B})$ of functions that are zero in $\partial\Omega$. For $\bar{v}\in P^0_0(\mathcal{B})$, the energy norm is defined by
\[
\Vert \bar{v}\Vert_E=\Vert G^{-1}\bar{v}\Vert_E.
\]

\bigskip

A basic result is the following,

\bigskip

\noindent\textbf{(2.2) Lemma. } Let $u,v\in P^1(\mathcal{T})$,  and $\bar{v}=g(v)\in P^0(\mathcal{B})$. Then 
\[
-\sum_{b\in\mathcal{B}}\int_{\partial b}\frac{\partial u}{\partial n}\bar{v}ds=\int_\Omega
\nabla u\cdot\nabla vdx,
\]
where $n$ is the outward pointing normal.

\noindent\textbf{Proof. }(Lemma 3 in\cite{BankRose})

\bigskip

An important observation is that the proof is local. That is, it suffices to show for $v=\phi_i$, $\bar{v}=\bar{\phi}_i$, namely
\[
\int_{\partial b_i}\frac{\partial u}{\partial n}ds=\int_{\Omega_i}
\nabla u\cdot\nabla \phi_idx.
\]

The equality is shown triangle by triangle for $t\in\mathcal{T}$ such that $t\subset\Omega_i$. 

Let us assume that the triangle $t=x_0x_1x_2$, it is proven that 
\[
-\int_{\partial b_0\cap t}\frac{\partial u}{\partial n}ds=\int_{\Omega_0\cap t}\nabla u\cdot\nabla\phi_0dx=
-d_2u_1+(d_1+d_2)u_0-d_1u_2.
\]
Here
\[
d_i=\frac{l_jl_kcos\, \theta_i}{4\vert t\vert},
\]
where $(i,j,k)$ is a cyclic permutation of $(0,1,2)$.

\subsection{The Poisson equation}

Consider the Dirichlet Boundary Value problem,
\[
-\Delta u(x)=f(x),\, x\in\Omega, \quad u(x)=0,\, x\in\partial\Omega,
\]
where $f\in L^2(\Omega)$.

Let $u$ be the $H^1_0(\Omega)$-weak solution of the Poisson equation. This solution satisfies
\[
a(u,v)=(f,v),\quad v\in H^1_0(\Omega).
\]

The FEM solution with linear elements is $u_L\in P^1_0(\mathcal{T})$ such that
\[
a(u_L,v)=(f,v),\quad v\in P^1_0(\mathcal{T}).
\] 
 
In light of the lemma, let us introduce the bilinear form
\[
\bar{a}(u,v)=-\sum_{b\in\mathcal{B}}\int_{\partial b}\frac{\partial u}{\partial n}\bar{v}ds.
\]

Define by $P^0_0(\mathcal{B})$ and $P^1_0(\mathcal{T})$ the subspaces of $P^0(\mathcal{B})$ and $P^1(\mathcal{T})$ whose elements are zero on $\partial\Omega$.

\bigskip

The box method consists on finding $u_B\in P^1_0(\mathcal{T})$ such that
\[
\bar{a}(u_B,\bar{v})=(f,\bar{v}),\quad \bar{v}\in P^0_0(\mathcal{B}).
\]

\bigskip

The method leads to the linear system
\[
\bar{A}\bar{U}=\bar{F}
\]
where, by the proposition
\[
\bar{A}_{ij}=\bar{a}(\phi_j,\bar{\phi}_i)=a(\phi_j,\phi_i),
\]
\[
\bar{F}_i=(f,\bar{\phi}_i)=\int_{b_i}f\, dx,
\]
\[
u_B(x)=\sum_i\bar{U}_i\phi_i(x).
\]

\bigskip

\noindent\textbf{(2.3) Remark. }The stiffness matrix is identical to that arising from the Finite Element Method with continuous piecewise linear approximations.
The corresponding linear systems differ only in the right hand side which are close in an average sense.

\section{FEM and BOX meshes and cell complexes}

In this section we describe the correspondence between FEM and BOX meshes and cell complexes.

\subsection{Primal Complex}

A $k-$dimensional simplex $\sigma^k\subset\mathbb{R}^d$ is the convex hull of $k+1$ geometrically independent points $x_0,x_1,\ldots,x_k\in \mathbb{R}^d$. 

Let $[x_0,x_1,\ldots,\hat{x}_i,\ldots x_k]$ the $(k-1)$-dimensional face with the $i$ vertex removed.

\bigskip

\noindent\textbf{(3.1) Definition. }A \emph{simplicial complex} $K\in\mathbb{R}^d$, is a collection of simplices in $\mathbb{R}^d$ such that
\begin{enumerate}
\item Every face of a simplex of $K$ is in $K$.

\item The intersection of any two simplices of $K$ is either a face of each of them or it is empty.
\end{enumerate}

The union of all simplices of $K$ treated as a subset of $\mathbb{R}^d$ is called the underlying space of $K$ and is denoted by $\vert K\vert$.

\bigskip

\noindent\textbf{(3.2) Definition. } A simplicial complex of dimension $n$ is called a \emph{manifold-like simplicial complex} if and only if $\vert K\vert$ is a $C^0-$manifold,  with or without boundary. More precisely
\begin{enumerate}
\item All simplices of dimension $k$ with $0\leq k\leq n-1$ must be a face of some simplex of dimension $n$ in $K$.

\item Each point on $\vert K\vert$ has a neighborhood homeomorphic to $\mathbb{R}^d$ or $n-$dimensional half-space.
\end{enumerate}

\bigskip

An oriented $k-$simplex is denoted by $\sigma^k=[x_0,x_1,\ldots, x_k]$. Its orientation corresponds to the orientation of the \emph{corner basis} at $x_0$. Namely, $\lbrace x_1-x_0,x_2-x_0,\ldots,x_k-x_0\rbrace$.  It is equivalent to state that two orderings of the vertices have the same orientation if they differ one from one another by an even permutation.

Let $[x_0,x_1,\ldots,\hat{x}_i,\ldots x_k]$ be a $(k-1)$-dimensional face. We say that the induced orientation on this face is the same of $\sigma^k$ if $i$ is even. If $i$ is odd, it is the opposite.

Assume that two oriented $k$-dimensional simplices, $\sigma$, $\tau$ share a face of dimension $k-1$. They have the same orientation if the induced orientation of the shared $k-1$ face induced by $\sigma$ is the opposite to that induced by  $\tau$. 

\bigskip

\noindent\textbf{(3.3) Definition. } A manifold-like simplicial complex $\Omega$ of dimension $n$ is called an \emph{oriented manifold-like
simplicial complex} if adjacent $n$-simplices (i.e., those that share a common $(n-1)$-face) have the same
orientation (orient the shared $(n-1)$-face oppositely) and simplices of dimensions $n-1$ and lower are
oriented individually. 

\bigskip

\noindent\textbf{2D FEM interpretation:}
\begin{itemize} 
\item The triangulation $\mathcal{T}$ may be regarded as a simplicial complex of dimension $2$ with underlying space $\vert\bar{\Omega}\vert$.

\item  For two adjacent triangles, orientation corresponds to set outward normals on the intersecting side,
opposite to each other.
\end{itemize}

\subsection{Dual Complex}

Given the triangulation $\mathcal{T}$, a point is selected in the closure of each triangle in the box method.  In the context of DEC, a natural choice is the circumcenter.

\bigskip

\noindent\textbf{(3.4) Definition. } The circumcenter of a $k$-simplex $\sigma^k$ is given by the center of the $k$-circumsphere. We will
denote the circumcenter of a simplex $\sigma^k$ by $c(\sigma^k)$. If the circumcenter of a simplex lies in its interior we call it a well-centered simplex. A simplicial complex all of whose simplices (of all dimensions) are well-centered will be called a well-centered simplicial complex.

\bigskip

Hereafter we shall consider triangulations corresponding to a well-centered simplicial complex.

\bigskip

\noindent\textbf{(3.5) Definition. } The circumcentric subdivision of a well-centered simplicial complex $K$ of dimension $n$
is denoted $csd K$, and it is a simplicial complex with the same underlying space as $K$ and consisting of all
simplices (each of which is called a subdivision simplex) of the form $[c(\sigma_1),\ldots,c(\sigma_k)]$ for $1\leq k\leq n$ (note that the index here is not dimension since it is a subscript). Here $\sigma_1\prec\sigma_2,\ldots,\prec\sigma_k$ 
(i.e.,$\sigma_i$ is a proper face  of $\sigma_j$ for all $i < j$) and the $\sigma_i$ are in $K$.

Each subdivision simplex in a given simplex $\sigma^k$ will be called a subdivision simplex of $\sigma^k$. Of these, a $j$ simplex ($j\leq k$) will be called a subdivision $j$-simplex of $\sigma^k$.

\bigskip

\noindent\textbf{(3.6) Circumcentric subdivision in 2D:} Consider a simplicial complex $K$ with vertices $x_0$, $x_1$ and $x_2$,
i.e., the complex consists of a triangle $x_0x_1x_2$, its edges and its vertices. Then $L=csd K$ consists of the
following elements:
\begin{itemize}
\item $L^{(0)}$ (the $0$-simplices of $csd K$): consists of the circumcenters $c(x_0) = x_0$, $c(x_1) = x_1$ and $c(x_2) = x_2$, the midpoints of the edges $c(x_0x_1)$, $c(x_1x_2)$ and $c(x_0x_2)$ and the circumcenter of the triangle, i.e.,
$c(v_0v_1v_2)$,

\item $L^{(1)}$ (the $1$-simplices of $csd K$): consists of $12$ edges, the two halves of each edge and edges joining
the circumcenter of the triangle to the vertices and midpoints of the edges,

\item $L^{(2)}$ (the $2$-simplices of $csd K$): consists of $6$ triangles, for instance

$x_0c(x_0x_1)c(x_0x_1x_2)$.
\end{itemize}

\bigskip

\noindent\textbf{(3.7) Definition. } Let $K$ be a well-centered manifold-like simplicial complex of dimension $n$ and let $\sigma^k$ be
one of its simplices. The circumcentric dual operator is given by 
\[
\star\sigma^k=\sum_{\sigma^k\prec\sigma^{k+1}\prec\ldots\prec\sigma^n}
\epsilon_{\sigma^k,\ldots,\sigma^n}[c(\sigma^k),c(\sigma^{k+1}),\ldots,\sigma^n],
\]
where the $\epsilon_{\sigma^k,\ldots,\sigma^n}$ coefficient ensure the the orientation of $[c(\sigma^k),c(\sigma^{k+1}),\ldots,\sigma^n]$ is
consistent with the orientation of the primal simplex, and the ambient volume-form.

\bigskip

The union of the interiors of the simplices in the definition of $\star\sigma^k$ is the (circumcentric) duall cell also denoted by $\star\sigma^k$. We will call each $(n-k)$-simplex $c(\sigma^k)c(\sigma{k+1})\ldots c(\sigma^n)$ an elementary dual simplex of $\sigma^k$. This is an $(n-k)$-simplex in $csd K$. The collection of dual cells is called the dual cell decomposition of $K$. This is dual cell complex and will be denoted $\star K$.

\bigskip

\noindent\textbf{2D BOX interpretation:} 

If $x_i$ is a vertex in the triangulation $\mathcal{T}$ of $\Omega$, then the dual cell of the $0-$simplex $[x_0]$
is the box $b_0$.

\section{DEC approximation of the Poisson equation}

We shall use the FEM-BOX notation when appropriate. 

\subsection{Chains, Forms, Integrals and Derivatives}

We start with an oriented manifold-like simplicial complex $K$ of dimension$n$.

\bigskip

\noindent\textbf{(4.1) Definition. }The space of $k-$chains of $k-$simplices is given by
\[
C_k(K)\equiv C_k(K,\mathbb{Z})=\left\lbrace
\sum_{j=1}^{N_k}\alpha_j\sigma^k_j:\alpha_j\in\mathbb{Z}
\right\rbrace,
\]
where $N_k$ is the number of $k-$simplices. $C_k(K,\mathbb{Z})$ is equipped with an Abelian group structure.

\bigskip

\noindent\textbf{(4.2) Definition. }The space of discrete $k-$forms (cochains) is the dual space of $k-$chains, namely
\[
C^k(K)=Hom(C_k(K),\mathbb{R}).
\]

\bigskip

\noindent\textbf{(4.3) Definition. }The integral of a $k-$cochain $\omega$ over a $k-$chain $c$ is defined by
\[
\int_c\omega^k=\omega^k(c).
\]

\bigskip

This bilinear pairing of cochain with chain is also denoted using bracket notation
\[
\langle\omega^k,c\rangle=\omega^k(c).
\]

The same construction can be applied to the dual cell complex$\star K$. We denote the Abelian group of $k-$chains on the dual complex by
\[
C_k(\star K)\equiv C_k(\star K,\mathbb{Z})=\left\lbrace
\sum_{j=1}^{N_k}\alpha_j\star\sigma^k_j:\alpha_j\in\mathbb{Z}
\right\rbrace.
\]
Also, the space of $k-$forms on the dual complex is denoted by
\[
C^k(\star K_d)=Hom(C_k(\star K),\mathbb{R}).
\]

\bigskip

\noindent\textbf{(4.4) Definition. }The $k$th boundary operator denoted by $\partial_k$ is the map
\[
\partial_k:C_k(K)\to C_{k-1}(K),
\]
given by
\[
\partial_k(\sigma^k)=\partial_k([x_0,x_1,\ldots,x_k])=\sum_{i=0}^k(-1)^i[x_0,\ldots,\hat{x}_i,\ldots,x_k],
\]
where $\hat{x}_i$ measn that $x_i$ is omitted.

\bigskip

\noindent\textbf{(4.5) Definition. }Let $\omega^k$ be a $k-$cochain. The $k$th discrete exterior derivative of $\omega$ is the transpose of the
$(k+1)-$st boundary operator
\[
\mathbb{D}_k=\partial_{k+1}^T.
\]

We have
\[
\langle\mathbb{D}_k\omega^k,\sigma^{k+1}\rangle=\langle\omega^k,\partial_{k+1}\sigma^{k+1}\rangle.
\]

\bigskip

The discrete exterior derivative operates on primal cochains. For dual cochains we have
\[
\mathbb{D}_{n-k}^{dual}=(-1)^k(\mathbb{D}_{k-1}^{primal})^T.
\]
The negative sign comes from the orientation on the dual mesh induced by the orientation on the primal mesh.
See \cite{Desbrun_K_T}.
\bigskip

Two more facts. It is readily seen that 
\[
\partial_k\circ\partial_{k+1}=0,
\]
and the Stokes' Theorem holds,
\[
\int_c\mathbb{D}_k\omega=\int_{\partial_{k+1}c}\omega.
\]

\subsection{The Discrete Hodge star operator and the Poisson Equation}

In the exterior calculus for smooth manifolds, the Hodge star, denoted $\ast$, is an isomorphism between the
space of $k-$forms and $(n-k)-$forms. It is apparent that a definition of a  discrete Hodge star ought to involve this relation in correspondence with the fact that the dual of a $k-$simplex is a $(n-k)-$cell. This is accomplished with the bilinear pairing.

\bigskip

\noindent\textbf{(4.6) Definition. }The discrete Hodge star operator is a map
\[
\ast:C^k(K)\to C^k(\star K).
\]

For a $k-$cochain $\omega^k$ it is defined as
\[
\frac{1}{\vert\star\sigma^k\vert}\langle\ast\omega^k,\star\sigma^k\rangle=
\frac{1}{\vert\sigma^k\vert}\langle\omega^k,\sigma^k\rangle
\]
for all $k-$simplices $\sigma^k$.

\bigskip

For all $k-$cochains $\omega^k$, it follows that
\[
\ast(\ast\omega^k)=(-1)^{k(n-k)}\omega^k.
\]

Also the analogue for $k-$chains is valid. Namely,
\[
\star(\star\sigma^k)=(-1)^{k(n-k)}\sigma^k.
\]

\bigskip

\noindent\textbf{The Poisson equation. }

\bigskip

Let us denote by $\vert A\vert$ the measure of a set $A$. Point sets are assigned unit measure. We restrict our discussion to the plane, 
and regard $\bar{\Omega}$ as a well-centered oriented manifold-like simplicial complex of dimension 2. 

Assume that $u$ and $f$ are $0-$primal forms. Consequently, a DEC approximation of the Poisson equation 
will result on nodal values for $u$.

Since $\ast f$ is a $2-$dual form we may write the Poisson equation in DEC formalism as follows,
\[
-\mathbb{D}_1^{dual}\ast\mathbb{D}_0^{primal}u=\ast f.
\]

Let $x_0$ be a vertex in the triangulation $\mathcal{T}$. It is also a $0-$simplex $[x_0]$. Applying the previous equation to its dual cell, we have
\[
\langle-\mathbb{D}_1^{dual}\ast\mathbb{D}_0^{primal}u,\star[x_0]\rangle=\langle\ast f,\star[x_0]\rangle.
\]

By the definition of Hodge star, the right hand side becomes
\[
\langle\ast f,\star[x_0]\rangle=\frac{\vert\star[x_0]\vert}{\vert [x_0] \vert}f_0.
\]
But $\vert [x_0] \vert=1$, thus
\[
\langle\ast f,\star[x_0]\rangle=\vert\star[x_0]\vert f_0.
\]

Notice that the dual cell of $[x_0]$  is exactly the box $b_0$,
\[
\star[x_0]=b_0.
\]

Hence, in box notation we have
\[
\langle\ast f,\star[x_0]\rangle=\vert b_0 \vert f_0.
\]

For the left hand side we prove the following

\bigskip

\noindent\textbf{(4.7) Theorem. }
\[
\langle-\mathbb{D}_1^{dual}\ast\mathbb{D}_0^{primal}u,\star[x_0]\rangle=
\sum_{[x_0,x_j,x_k]\subset\Omega_0}\left[
-\frac{e_k}{l_k}u_j+\left(\frac{e_j}{l_j}+ \frac{e_k}{l_k}\right)u_0
-\frac{e_j}{l_j}u_k
\right].
\]
\noindent\textbf{Proof. }

Discrete derivatives are the transpose of boundary operators, thus
\[
\langle-\mathbb{D}_1^{dual}\ast\mathbb{D}_0^{primal}u,\star[x_0]\rangle=
-\langle\ast\mathbb{D}_0^{primal}u,\partial_2^{dual}\star[x_0]\rangle.
\]

It is apparent that 
\[
\partial_2^{dual}\star[x_0]=\eth b_0=e_1\cup e_2.
\]

By linearity of the boundary operator and the integral, it is enough to show the equality for a triangle $t\in\mathcal{T}$, $t\subset\Omega_0$. 
Let $t_0=x_0x_1x_2$ be one of such triangles, or the $2-$primal simplex $[x_0,x_1,x_2]$. 

Let us discuss the orientation of the dual $1-$simplices $[c([x_0,x_1]),c([x_0,x_1,x_2])\equiv [m_2,C]$ and
$[c([x_0,x_2]),c([x_0,x_1,x_2])\equiv [m_1,C]$. For the former, consider the $2-$simplex $[x_0,m_2,C]$. This is related to the area form
$dx^1\wedge dx^2$, up to a sign determined by the relative orientation of $[x_0,m_2]$ and $[x_0,x_1]$. Thus we have that
\[
dx^1\wedge dx^2\sim sgn([x_0,m_2],[x_0,x_1])[x_0,m_2,C]
\]
Thus we have that the correct orientation for the simplex $[m_2,C]$ is given by
\[
sgn([x_0,m_2],[x_0,x_1])\cdot sgn([x_0,m_2,C],[x_0,x_1,x_2])\equiv (+1)(+1)=+1.
\]

Similarly,
\[
dx^1\wedge dx^2\sim sgn([x_0,m_1],[x_0,x_2])[x_0,m_1,C]
\]
and the orientation for the simplex $[m_1,C]$ is given by
\[
sgn([x_0,m_1],[x_0,x_2])\cdot sgn([x_0,m_1,C],[x_0,x_1,x_2])\equiv (+1)(-1)=-1.
\]

We are led to compute in the simplex $[x_0,x_1,x_2]$
\[
-\langle\ast\mathbb{D}_0^{primal}u,[m_2,C]-[m_1,C]\rangle= 
-\langle\ast\mathbb{D}_0^{primal}u,[m_2,C]\rangle
-\langle\ast\mathbb{D}_0^{primal}u,-[m_1,C]\rangle
\]

By the definition of the dual cell operator
\[
-\langle\ast\mathbb{D}_0^{primal}u,[m_2,C]-[m_1,C]\rangle= 
-\langle\ast\mathbb{D}_0^{primal}u,\star[x_0,x_1]\rangle
-\langle\ast\mathbb{D}_0^{primal}u,-\star[x_2,x_0]\rangle
\]
Applying the definition of the discrete Hodge star operator we obtain
\[
\begin{array}{rcl}
-\langle\ast\mathbb{D}_0^{primal}u,[m_2,C]-[m_1,C]\rangle & = &
-\frac{\vert \star[x_0,x_1]\vert}{\vert [x_0,x_1]\vert}\langle\mathbb{D}_0^{primal}u,[x_0,x_1]\rangle \\
 & & \\
 & & +\frac{\vert \star[x_2,x_0]\vert}{\vert [x_2,x_0]\vert}\langle\mathbb{D}_0^{primal}u,[x_2,x_0]\rangle \\
 & & \\
  & \equiv &
  -\frac{e_2}{l_2}\langle\mathbb{D}_0^{primal}u,[x_0,x_1]\rangle
+\frac{e_1}{l_1}\langle\mathbb{D}_0^{primal}u,[x_2,x_0]\rangle
\end{array}
\]

Now 
\[
\begin{array}{rcl}
-\langle\ast\mathbb{D}_0^{primal}u,[m_2,C]-[m_1,C]\rangle & = &
  -\frac{e_2}{l_2}\langle u,\partial_1^{primal}[x_0,x_1]\rangle
+\frac{e_1}{l_1}\langle u,\partial_1^{primal}[x_2,x_0]\rangle \\
& & \\
  & = &
  -\frac{e_2}{l_2}(u_1-u_0)
+\frac{e_1}{l_1}(u_0-u_2) \\
& & \\
  & = &
 -\frac{e_2}{l_2}u_1+\left(\frac{e_1}{l_1}+ \frac{e_2}{l_2}\right)u_0
-\frac{e_1}{l_1}u_2.
\end{array}
\]
The result follows.$\blacksquare$

\bigskip

Consequently, in the dual cell $\star[x_0]\equiv b_0$ we have
\[
\sum_{[x_0,x_j,x_k]\subset\Omega_0}\left[
-\frac{e_k}{l_k}u_j+\left(\frac{e_j}{l_j}+ \frac{e_k}{l_k}\right)u_0
-\frac{e_j}{l_j}u_k
\right]
=\vert b_0 \vert f_0
\]

Assembling these block submatrices we obtain $u_D$, the DEC solution of this system with zero boundary conditions.

\subsection{Error estimates  in energy norm}

We show that the stiffness matrix of de BOX and DEC methods coincide.

\bigskip

Let us consider the triangle in Figure \ref{fig:example}. Using the fact that

\begin{figure}[htbp] %  figure placement: here, top, bottom, or page
   \centering
   \includegraphics[width=2in]{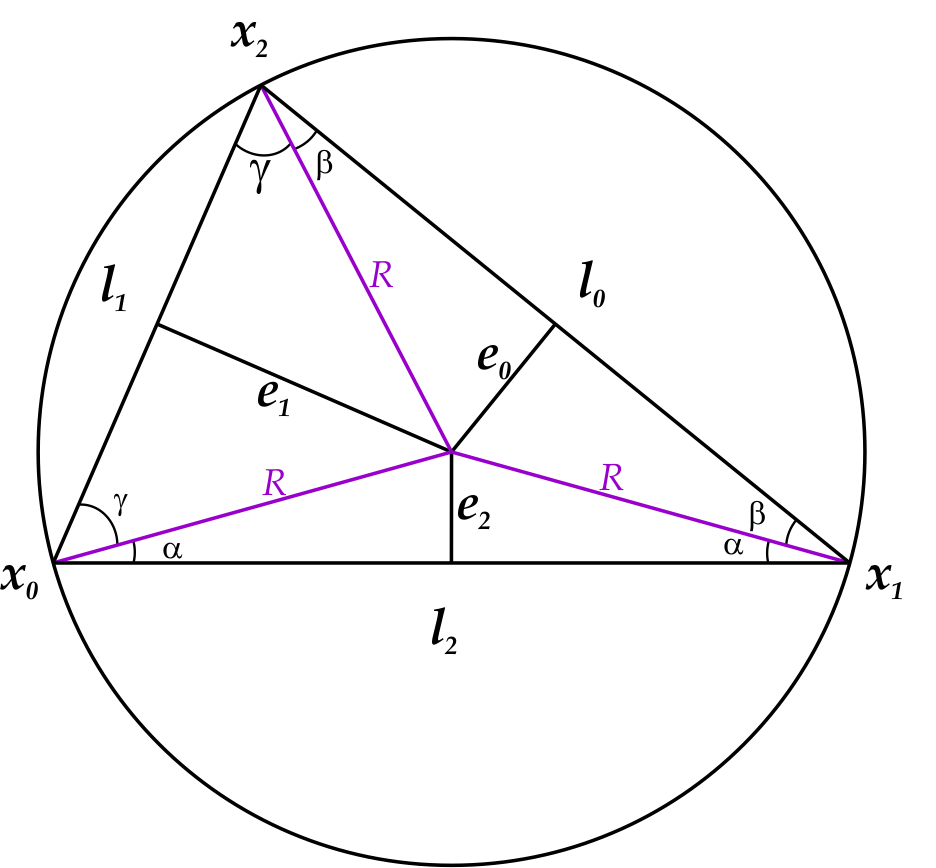} 
  \caption{}
   \label{fig:example}
\end{figure}

\[
d_0=\frac{l_1l_2cos\, \theta_0}{4\vert t\vert},
\]
and
\[
\vert t \vert=\frac{l_0l_1l_3}{4R},
\]
we obtain by elementary geometry that 
\begin{equation}
d_0=\frac{e_0}{l_0}.
\label{d_e_l}
\end{equation}

It is apparent that this is valid for every vertex in $\mathcal{T}.$ We are led to the following result.

\bigskip

\noindent\textbf{(4.8) Theorem. }Let $f\in L^2(\Omega)$ be piecewise constant with respecto to $\mathcal{B}$.
Then, there is a positive constant $M$ such that
\[
\Vert u-u_L\Vert_E\leq \Vert u-u_D\Vert_E\leq M\Vert u-u_L\Vert_E.
\]
\noindent\textbf{Proof. }

Multiplying the Poisson equation by $v\in\mathcal{S}^0(\mathcal{B})$ and integrating by parts we obtain
\[
\sum_{b\in\mathcal{B}}\int_b\nabla u\cdot\nabla vdx +\bar{a}(u,v)=(f,v).
\]
Recall that the box solution $u_B\in P^1_0(\mathcal{T})$ satisfies
\[
\bar{a}(u_B,\bar{v})=(f,\bar{v})\,\quad \bar{v}\in P^0_0(\mathcal{B}).
\]

Thus , for $f\in L^2(\Omega)$ piecewise constant with respecto to $\mathcal{B}$ we have from equation (\ref{d_e_l}) that $u_b\equiv u_D$, so $u_D$ satisfies
\[
\bar{a}(u_D,\bar{v})=(f,\bar{v})\,\quad \bar{v}\in P^0_0(\mathcal{B}).
\]
It follows that
\[
\bar{a}(u-u_D,\bar{v})=0.
\]
Finally, continue with the proof of Theorem 1 and Corollary 1 in \cite{BankRose}.$\blacksquare$

\subsection{General triangulations}

Let us assume that the triangle $t=x_0x_1x_2$. In the box method, the line integral 
\[
\int_{\eth b_0\cap t}\frac{\partial u}{\partial n}ds
\]
depends only on $m_2$ and $m_1$, the endpoints of the integration path. So the location of the distinguished point $C\in\bar{t}$ is arbitrary. 

For the DEC method, we argue that the circumcenter need not be in the interior of $t$.

Starting with the FEM local matrix computation, it is straightforward to see that

\begin{eqnarray*}
{e_0\over l_0}
 &=&
{1\over l_0^2}
\det\left(\begin{array}{ccc}
x_1 & y_1 & 1\\
x_c & y_c & 1\\
x_2 & y_2 & 1
                            \end{array}
\right),\\
{e_1\over l_1}
 &=&
{1\over l_1^2}
\det\left(\begin{array}{ccc}
x_2 & y_2 & 1\\
x_c & y_c & 1\\
x_0 & y_0 & 1
                            \end{array}
\right).\\
{e_2\over l_2}
 &=&
{1\over l_2^2}
\det\left(\begin{array}{ccc}
x_0 & y_0 & 1\\
x_c & y_c & 1\\
x_1 & y_1 & 1
                            \end{array}
\right),\\
\end{eqnarray*}

These expressions are valid regardless of the location of the circumcenter and can, indeed, take negative values.
Such signs are essential for the calculations to work for general meshes. In  \cite{Herrera_etal} we 
tested the  DEC method for bad quality meshes, in the FEM jargon. As expected, The numerical performance 
was of the first order in these not well centered meshes.

\end{document}